\documentclass[11pt, a4paper]{amsart}
\usepackage{amsmath, amsthm, amssymb}
\usepackage{txfonts, mathrsfs}
\usepackage{graphicx}
\usepackage{setspace}
\usepackage[left=3.8cm, right=3.8cm, top=3.4cm, bottom=3.2cm]{geometry}

\newtheorem*{thm*}{Theorem}

\newtheorem{thm}{Theorem}[section]
\newcommand{\bt}{\begin{thm}}
\newcommand{\et}{\end{thm}}

\newtheorem{cor}[thm]{Corollary}
\newcommand{\bc}{\begin{cor}}
\newcommand{\ec}{\end{cor}}

\newtheorem{lem}[thm]{Lemma}
\newcommand{\bl}{\begin{lem}}
\newcommand{\el}{\end{lem}}

\newtheorem{prop}[thm]{Proposition}
\newcommand{\bp}{\begin{prop}}
\newcommand{\ep}{\end{prop}}

\newtheorem{defn}[thm]{Definition}
\newcommand{\bd}{\begin{defn}}      
\newcommand{\ed}{\end{defn}}

\newtheorem{rmrk}[thm]{Remark}
\newcommand{\br}{\begin{rmrk}}
\newcommand{\er}{\end{rmrk}}

\newtheorem{quest}[thm]{Question}
\newcommand{\bq}{\begin{quest}}
\newcommand{\eq}{\end{quest}}

\newcommand{\N}{\mathbb{N}}

\newcommand{\R}{\mathbb{R}}

\newdimen\vintkern\vintkern12pt

\def\vint{-\kern-\vintkern\int}

\newcommand{\hm}{{\mathcal H}}
\newcommand{\dist}{\operatorname{dist}}
\newcommand{\diam}{\operatorname{diam}}
\newcommand{\trace}{\operatorname{tr}}
\newcommand{\Area}{\operatorname{Area}}
\newcommand{\osc}{\operatorname{osc}}

\newcommand{\Image}{\operatorname{im}}

\newcommand{\jac}{{\mathbf J}}
\newcommand{\ap}{\operatorname{ap}}
\newcommand{\apmd}{\operatorname{md}}

\DeclareMathOperator{\interior}{int}
\newcommand{\id}{\operatorname{id}}

\newcommand{\rel}{\text{\;rel\;}}

\DeclareMathOperator{\MOD}{mod}

\begin{document}
\bibliographystyle{plain}

\title[Quasiconformal uniformization of metric surfaces]{Quasiconformal uniformization of metric surfaces of higher topology}

\author{Damaris Meier}


\address
  {Department of Mathematics\\ University of Fribourg\\ Chemin du Mus\'ee 23\\ 1700 Fribourg, Switzerland}
\email{damaris.meier@unifr.ch}

\date{\today}

\begin{abstract}
    We establish the following uniformization result for metric spaces $X$ of finite Hausdorff 2-measure. If $X$ is homeomorphic to a smooth 2-manifold $M$ with non-empty boundary, then we show that $X$ admits a quasiconformal almost parametrization $M\to X$, by only assuming that $X$ is locally geodesic and has rectifiable boundary. In particular, we recover a corollary of Ntalampekos and Romney by using the solution of the Plateau problem. After putting additional assumptions on $X$, we show that the quasiconformal almost parametrization upgrades to a quasisymmetry or a geometrically quasiconformal map, implying statements analogous to the uniformization theorems of Bonk and Kleiner as well as Rajala for surfaces of higher topology.
\end{abstract}

\maketitle

\section{Introduction}
\subsection{Background and statement of main result}
The classical uniformization theorem states that any Riemann surface is conformally equivalent to a surface of constant curvature. One of the main questions in the field of analysis on metric spaces asks under which conditions on a metric space $X$ homeomorphic to some model space $M$ there exists a homeomorphism $u\colon M\to X$ with good geometric and analytic properties.

The first breakthrough result regarding uniformization of metric surfaces is due to Bonk and Kleiner \cite{BK02} and asserts that if a metric space $X$ is Ahlfors $2$-regular and homeomorphic to $S^2$, then $X$ is quasisymmetrically equivalent to $S^2$ if and only if $X$ is linearly locally connected.
For the definitions we refer to Section~\ref{sec:applications}. In this work, a \emph{smooth surface} refers to a smooth compact oriented and connected Riemannian $2$-dimensional manifold with possibly non-empty boundary and a \emph{metric surface} is a metric space homeomorphic to a smooth surface.

Let $X$ be a metric space of locally finite Hausdorff $2$-measure. The \emph{modulus of a curve family $\Gamma$} in $X$ is defined by $$\MOD(\Gamma):= \inf_\rho \int_X\rho^2\,d\hm^2,$$ where the infimum is taken over all Borel functions $\rho\colon X\to [0,\infty]$ for which $\int_\gamma\rho\geq 1$ holds for every $\gamma\in\Gamma$, see Section~\ref{sec:prelimsMod}. A homeomorphism $u\colon M\to X$ is \emph{geometrically quasiconformal} if it leaves the modulus of curve families quasi-invariant. Rajala \cite{Raj14} showed that every metric space $X$ homeomorphic to $\R^2$ of locally finite Hausdorff $2$-measure admits a geometrically quasiconformal map $u$ from a domain $\Omega\subset\R^2$ to $X$ if and only if $X$ satisfies a condition called reciprocality, see Section~\ref{sec:applications}. 

In a next step, Ntalampekos and Romney \cite{NR21} as well as Wenger and the author \cite{MW21}, showed independently and with two different approaches the following result. If $X$ is locally geodesic and one relaxes the assumptions on $u$, then the condition of reciprocality can be dropped completely. The assumptions on $u$ are such that $u$ is a continuous, monotone surjection satisfying  
    \begin{align}
        \label{ineq:modulus}
        \MOD(\Gamma)\leq K\cdot \MOD(u\circ\Gamma)
    \end{align}
for $K\geq 1$ and every family $\Gamma$ of curves in $M$, where we denote by $u\circ\Gamma$ the family of curves $u\circ\gamma$ with $\gamma\in\Gamma$. Here, a map $u\colon M\to X$ is \emph{monotone} if the preimage of a point is connected and equivalently, if $u$ is a uniform limit of homeomorphisms $M\to X$ (see \cite{You48} and Proposition~\ref{prop:genOf2.9} below). Whenever $u\colon M\to X$ is in addition a homeomorphism, then by \cite{Wil12}, condition \eqref{ineq:modulus} is equivalent to the so-called analytic definition of quasiconformality. A similar relation holds for $u$ being continuous, monotone and satisfying \eqref{ineq:modulus}, see \cite[Theorem~7.1]{NR21}.

All results mentioned so far were stated for simply connected surfaces. It is very natural to wonder to what extent these statements can be extended to surfaces of higher topology. Rajala's uniformization theorem has been generalized in \cite{Iko19}, while generalizations of the theorem of Bonk-Kleiner can be found in \cite{MW13}, \cite{GW18} and \cite{FM21}. 

Very recently, Ntalampekos and Romney \cite{NR22} showed that the assumption of being locally geodesic in \cite{NR21} and \cite{MW21} is not needed. In particular, their result holds for all metric surfaces $X$ of locally finite Hausdorff 2-measure, even for surfaces of higher genus and with non-empty boundary. 

In this work we focus on the approach from \cite{MW21}, which is closely related to the existence and regularity of energy and area minimizing discs in metric spaces admitting a quadratic isoperimetric inequality developed by Lytchak and Wenger in \cite{LW15-Plateau}. Moreover, we will make use of the theory of area minimizing surfaces in homotopy classes in metric spaces \cite{SW20}. Our main result is the following Corollary of \cite[Theorem~1.3]{NR22}, using a completely different proof strategy than in \cite{NR22}.

\bt\label{thm:main}
    Let $X$ be a locally geodesic metric space homeomorphic to a smooth surface $M$ with non-empty boundary. If $X$ is of finite Hausdorff $2$-measure and has rectifiable boundary, then there exists a Riemannian metric $g$ on $M$ and a continuous, monotone surjection $u\colon M\to X$ such that
    \begin{align*}
        \MOD(\Gamma)\leq K\cdot \MOD(u\circ\Gamma)
    \end{align*}
    holds for every family $\Gamma$ of curves in $M$ with $K=\frac{4}{\pi}$.
\et

The Riemannian metric $g$ can be chosen in such a way that it is of constant sectional curvature $-1$, $0$ or $1$ and $\partial X$ is geodesic. Note that the constant $\frac{4}{\pi}$ is optimal (see \cite[Theorem~1.3]{Iko19}) and that in this generality, there doesn't have to exist a homeomorphism satisfying \eqref{ineq:modulus}, see e.g. \cite[Example~3.2]{MW21}.

\subsection{Proof of main theorem} 

Let $X$ be a metric space homeomorphic to a smooth surface $M$. The boundary $\partial X$ is the disjoint union of $k\geq0$ Jordan curves in $X$. Denote by $[\partial X]$ the set of all weakly monotone parametrizations of $\partial X$, i.e. uniform limits of homeomorphisms from $S$ to $\partial X$, where $S$ is any space homeomorphic to the disjoint union of $k$ copies of $S^1$. Let $\Lambda(M,\partial X,X)$ be the possibly empty family of Sobolev maps $u\in N^{1,2}(M,X)$ such that the trace $\trace(u)$ has a continuous representative in $[\partial X]$. Consider a homeomorphism $\varphi\colon M\to X$ and let $h\colon K\to M$ be a triangulation of $M$, see Section~\ref{sec:prelimsRelHom}. Denote by $K^1$ the $1$-skeleton of $K$ and by $\partial K\subset K^1$ the preimage of $\partial M$ under $h$. Two continuous maps $\varrho,\varrho'\colon K^1\to X$ with $\varrho|_{\partial K},\varrho'|_{\partial K}\in[\partial X]$ are said to be \emph{homotopic relative to $\partial X$} if there exists a homotopy $H$ between $\varrho$ and $\varrho'$ with $H(\cdot,t)|_{\partial K}\in[\partial X]$ for all $t\in[0,1]$. The common relative homotopy class is denoted by $[\varrho]_{\partial X}$. Note that if $\partial X$ is empty, then $[\varrho]_{\partial X}$ corresponds to the usual homotopy class of $\varrho$. If $u\in \Lambda(M,\partial X,X)$ is continuous, we say that \emph{$u$ is $1$-homotopic to $\varphi$ relative to $\partial X$}, denoted $u\sim_1\varphi\rel\partial X$, if $$[u\circ h|_{K^1}]_{\partial X}=[\varphi\circ h|_{K^1}]_{\partial X}$$ for some and thus any triangulation $h$ on $M$, see \cite{SW20}. The $1$-homotopy class $u_{\#,1}[h]$ of an arbitrary $u\in\Lambda(M,\partial X,X)$ will be defined in Section~\ref{sec:prelimsRelHom} using the existence of a local quadratic isoperimetric inequality in the $\varepsilon$-thickening $X_{\varepsilon}$ of $X$ and a suitable retraction $X_{\varepsilon}\to X$. In a first step we show that the family
$$\Lambda(M,\varphi,X):=\{u\in\Lambda(M,\partial X,X):u\sim_1\varphi\rel\partial X\}$$
is not empty. Notice that in \cite{SW20} the existence of a map in $\Lambda(M,\varphi,X)$ highly depends on the fact that $X$ admits a local quadratic isoperimetric inequality. In this article we make use of the 2-dimensional structure of $X$ to prove the following generalization of \cite[Theorem~1.4]{MW21}.
\bt\label{thm:existence-Sobolev}
    Let $X$ be a locally geodesic metric space homeomorphic to a smooth surface $M$ that is not a sphere and let $\varphi\colon M\to X$ be a homeomorphism. If $\hm^2(X)<\infty$ and $\ell(\partial X)<\infty$ then the family $\Lambda(M,\varphi,X)$ is not empty.
\et

For convenience, we denote by 
$\Lambda_{\text{metr}}(M,\varphi,X)$ the family of pairs $(u,g)$, where $u\in\Lambda(M,\varphi,X)$ and $g$ is a Riemannian metric on $M$. Moreover, a pair $(u,g)\in\Lambda_{\text{metr}}(M,\varphi,X)$ satisfying 
    \begin{align*}
        E_+^2(u,g)=\inf\{E_+^2(v,h):v\in\Lambda(M,\varphi,X),\,h\text{ is a Riemannian metric on }M\},
    \end{align*} 
is called \emph{energy minimizing}. Here, $E_+^2$ denotes the Reshetnyak energy, see Section \ref{sec:prelimsSobolev}. Assuming that $X$ is locally geodesic, has rectifiable boundary and is of finite Hausdorff $2$-measure, then by Theorem~\ref{thm:existence-Sobolev}, the family $\Lambda(M,\varphi,X)$ is not empty. Theorem~\ref{thm:Sobolev-energy-min} below implies the existence of an energy minimizing pair $(u,g)\in\Lambda_{\text{metr}}(M,\varphi,X)$. The regularity of such an energy minimizer follows from the next theorem generalizing \cite[Theorem~1.3]{MW21}.

\bt\label{thm:continuity-of-energy-min-Sobolev}
    Let $M$ be a smooth surface with non-empty boundary, $X$ a locally geodesic metric space homeomorphic to $M$ and $\varphi\colon M\to X$ a homeomorphism. If $(u,g)\in\Lambda_{\text{metr}}(M,\varphi,X)$ is an energy minimizing pair, then $u$ has a representative which is continuous and extends continuously to the boundary.
\et

Hence, we obtain the existence of an energy minimizing pair $(u,g)\in\Lambda_{\text{metr}}(M,\varphi,X)$ with $u$ being continuous. In a next step we show that every such $u$ is a uniform limit of homeomorphisms and therefore monotone, compare to \cite[Theorem~1.2]{LW20-param}.

\bt\label{thm:cont-and-inf-iso-implies-monotone}
    Let $X$ be a locally geodesic metric space homeomorphic to a smooth surface $M$ with non-empty boundary and let $\varphi\colon M\to X$ be a homeomorphism. If a continuous map $u\in\Lambda(M,\varphi, X)$ and a Riemannian metric $g$ on $M$ satisfy 
    \begin{align*}
        E_+^2(u,g)=\inf\{E_+^2(v,h):v\in\Lambda(M,\varphi,X),\,h\text{ is a Riemannian metric on }M\},
    \end{align*} 
    then $u$ is a uniform limit of homeomorphisms from $M$ to $X$.
\et

By \cite[Corollary~1.3]{FW20}, the map $u$ is infinitesimally isotropic with respect to $g$ and thus infinitesimally $K$-quasiconformal with respect to $g$ for $K=\frac{4}{\pi}$ (see \cite{LW16-energy-area}). Consider Section~\ref{sec:prelimsSobolev} for the definitions of infinitesimal isotropy and infinitesimal quasiconformality. By arguing exactly as in the proof of \cite[Corollary~3.5]{LW20-param}, we obtain that $u$ satisfies (\ref{ineq:modulus}) with $K=\frac{4}{\pi}$. This establishes Theorem~\ref{thm:main}.\\

The paper is structured as follows. In Section~\ref{sec:prelims} we assemble necessary definitions and results needed later on. Theorem~\ref{thm:existence-Sobolev} is shown in Section~\ref{sec:existence}, while Sections~\ref{sec:continuity} and \ref{sec:monotone} are devoted to the proofs of Theorem~\ref{thm:continuity-of-energy-min-Sobolev} and Theorem~\ref{thm:cont-and-inf-iso-implies-monotone}, respectively. In Section~\ref{sec:applications} we show that the map $u$ from Theorem~\ref{thm:main} upgrades to a quasisymmetric or a geometrically quasiconformal map after further assumptions on the underlying spaces, recovering the above mentioned generalizations of Bonk-Kleiner and Rajala.  \\

{\bf Acknowledgments:} I wish to thank my PhD advisor Stefan Wenger for many important inputs and interesting discussions regarding this work.

\section{Preliminaries}\label{sec:prelims}

\subsection{Basic definitions and notations}

Let $(X,d)$ be a metric space. We denote the \emph{open ball} in $X$ of radius $r>0$ centered at a point $x\in X$ by $B(x,r)$. The \emph{open and closed unit discs} in $\R^2$ are given by
$$D:=\{z\in\R^2:|z|<1\},\qquad \overline{D}:=\{z\in\R^2:|z|\leq1\},$$
where $|\cdot|$ is the Euclidean norm. A set $\Omega\subset X$ homeomorphic to the unit disc $D$ is a \emph{Jordan domain} in $X$ if its boundary $\partial\Omega\subset X$ is a \emph{Jordan curve} in $X$, i.e. a subset of $X$ homeomorphic to $S^1$. The \emph{length} of a curve $c$ in $X$ is denoted by $\ell(c)$. A curve $c$ is said to be \emph{rectifiable} if $\ell(c)<\infty$ and \emph{locally rectifiable} if each of its compact subcurves is rectifiable. Moreover, a curve $c\colon[a,b]\to X$ is called \emph{geodesic} if $\ell(c)=d(c(a),c(b))$. A metric space $(X,d)$ is \emph{geodesic} if every pair of points in $X$ can be joined by a geodesic in $X$ and it is called \emph{locally geodesic} if every point $x\in X$ has a neighborhood $U$ such that any two points in $U$ can be joined by a geodesic in $X$.
 
For $s\geq 0$, we denote the \emph{$s$-dimensional Hausdorff
measure} of a set $A\subset X$ by $\hm^s(A)$. The normalizing constant is chosen in such a way that if $X$ is the Euclidean space $\R^n$, the Lebesgue measure agrees with $\hm^n$ on open subsets of $\R^n$. If $(M,g)$ is a Riemannian manifold of dimension $n$ then the $n$-dimensional Hausdorff measure $\hm_g^n$ on $(M,g)$ coincides with the Riemannian volume. We emphasize that throughout this paper, the reference measure on metric spaces will always be the $2$-dimensional Hausdorff measure. 

Let $g$ be a smooth Riemannian metric on a smooth surface $M$ such that the boundary of $M$ is geodesic with respect to $g$. We call the metric $g$ \emph{hyperbolic} if it is of constant sectional curvature $-1$, and \emph{flat} if it has vanishing sectional curvature as well as an associated Riemannian 2-volume satisfying $\hm^2_g(M)=1$.

\subsection{Conformal modulus}
\label{sec:prelimsMod}
Let $X$ be a metric space of locally finite Hausdorff $2$-measure and $\Gamma$ a family of curves in $X$. A Borel function $\rho\colon X\to [0,\infty]$ is said to be \emph{admissible for $\Gamma$} if $\int_\gamma \rho\geq 1$ for every locally rectifiable curve $\gamma\in\Gamma$. For the definition of the path integral $\int_\gamma\rho$ see \cite{HKST15}. The \emph{modulus} of the curve family $\Gamma$ is now defined by $$\MOD(\Gamma):= \inf_\rho\int_X\rho^2\,d\hm^2,$$ where the infimum is taken over all admissible functions for $\Gamma$. If $\Gamma$ contains a constant curve, then $\MOD(\Gamma)=\infty$. We say that a property holds \emph{for almost every curve in $\Gamma$} if it holds for every curve in $\Gamma_0$ for some $\Gamma_0\subset \Gamma$ with $\MOD(\Gamma\setminus \Gamma_0)=0$.

A homeomorphism $u\colon M\to X$ is \emph{geometrically quasiconformal} if there exists $K\geq 1$ such that 
$$ K^{-1}\cdot \MOD(\Gamma)\leq \MOD(u\circ\Gamma)\leq K\cdot \MOD(\Gamma)$$
for every family $\Gamma$ of curves in $M$.

\subsection{Metric space valued Sobolev maps}
\label{sec:prelimsSobolev}
We now give a brief overview over some basic concepts used in the theory of metric space valued Sobolev maps based on upper gradients. Note that several other equivalent definitions of Sobolev spaces exist. For more details consider e.g. \cite{HKST15}. 

Let $(X,d)$ be a complete metric space and $M$ a smooth surface. Fix a Riemannian metric $g$ on $M$ and consider a domain $U\subset M$. Let $u\colon U\to X$ be a map and $\rho\colon U\to [0,\infty]$ a Borel function. Then, $\rho$ is called \emph{(weak) upper gradient of $u$ with respect to $g$} if
\begin{align}
    \label{ineq:UpperGradient}
    d(u(\gamma(a)),u(\gamma(b)))\leq\int_{\gamma}\rho(s)\;ds
\end{align}
for (almost) every rectifiable curve $\gamma\colon[a,b]\to U$.

Denote by $L^2( U,X)$ the family of measurable essentially separably valued maps $u\colon U\to X$ such that the function $u_x(z):=d(u(z),x)$ is in the space $L^2( U)$ of $2$-integrable functions for some and hence any $x\in X$. A sequence $(u_k)\subset L^2( U,X)$ is said to \emph{converge in $L^2( U,X)$} to a map $u\in L^2( U,X)$ if
$$\int_{ U}d^2(u_k(z),u(z))\;d\hm^2_g(z)\to0$$
as $k$ tends to infinity. The \emph{(Newton-)Sobolev space} $N^{1,2}( U,X)$ is the collection of maps $u\in L^2( U,X)$ such that $u$ has a weak upper gradient in $L^2( U)$. Every such $u$ has a minimal weak upper gradient denoted by $\rho_u$, meaning that $\rho_u\in L^2( U)$ and for every weak upper gradient $\rho$ of $u$ in $L^2( U)$ it holds that $\rho_u\leq\rho$ almost everywhere on $ U$. Moreover, $\rho_u$ is unique up to sets of measure zero (see e.g. \cite[Theorem~6.3.20]{HKST15}). We emphasize that the definition of $N^{1,2}( U,X)$ is independant of the chosen metric $g$ on $M$.

The \emph{Reshetnyak energy} of a map $u\in N^{1,2}( U,X)$ with respect to $g$ is defined by
$$E^2_+(u,g):=\int_{ U}\rho_u(z)^2\;d\hm_g^2(z).$$
Note that this definition of energy agrees with the one given in \cite[Definition 2.2]{FW21}; in particular, $E_+^2$ is invariant under precompositions with conformal diffeomorphisms.\\

Consider a domain $V\subset\R^2$. A map $v\colon V\to X$ is said to be \emph{approximately metrically differentiable at $z\in V$} if there is a necessarily unique seminorm $s$ on $\R^2$ such that
$$\ap\lim_{y\to z}\frac{d(v(y),v(z))-s(y-z)}{|y-z|}=0,$$
where $\ap\lim$ denotes the approximate limit (see e.g. \cite{EG92}). If such a seminorm exists, it is called \emph{approximate metric derivative of $v$ at $z$}, denoted $\apmd v_z$. Consider an open set $W\subset\R^2$, a point $w\in W$ and a diffeomorphism $\varphi\colon W\to V$. If the map $v\colon V\to X$ is approximately metrically differentiable at $\varphi(w)$ then the composition $v\circ\varphi$ is approximately metrically differentiable at $w$ with
$$\apmd(v\circ\varphi)_w=\apmd v_{\varphi(w)}\circ d\varphi_w.$$
Along with \cite[Proposition~4.3]{LW15-Plateau} this implies that if $u\in N^{1,2}( U,X)$ then for almost every $z\in U $ the composition $u\circ\psi^{-1}$ is approximately metrically differentiable at $\psi(z)$ for an arbitrary chart $(U,\psi)$ around $z$. Define the seminorm $\apmd u_z$ on $T_z M$ by
$$\apmd u_z:=\apmd(u\circ\psi^{-1})_{\psi(z)}\circ d\psi_z.$$
Note that this definition is independent of the choice of chart and $\apmd u_z$ is called \emph{approximate metric derivative of $u$ at $z$}.

The \emph{jacobian} $\jac(s)$ of a seminorm $s$ on $\R^2$ is defined to be the Hausdorff $2$-measure on $(\R^2,s)$ of the unit square if $s$ is a norm and zero otherwise. By identifying $(T_zM,g(z))$ with $(\R^2,|\cdot|)$ via a linear isometry, we are able to define the jacobian of a seminorm $s$ on $T_zM$.

\bd
    The \emph{parametrized (Hausdorff) area} of $u\in N^{1,2}( U,X)$ is given by
    $$\Area(u):=\int_{ U}\jac(\apmd u_z)\;d\hm^2_g(z).$$
\ed

We emphasize that the parametrized area of $u\in N^{1,2}( U,X)$ is invariant under precompositions with biLipschitz homeomorphisms, and thus independent of the Riemannian metric $g$. Moreover, if $u$ is a homeomorphism onto its image then the jacobian $\jac(\apmd u_z)$ agrees with the Radon-Nikodym derivative of the measure $u^*\hm^2(B):= \hm^2(u(B))$ with respect to the Lebesgue measure at almost every point $z\in U$.

Recall that by John's theorem (see e.g. \cite[Theorem~2.18]{AlvT04}), the unit ball $B$ of a $2$-dimensional space equipped with a norm $s$ contains a unique ellipse $E$ of maximal area, called \emph{John's ellipse of $s$}. The $\mu^i$-jacobian $\jac_{\mu^i}(s)$ of a semi-norm $s$ on $\R^2$ is given by $\jac_{\mu^i}(s)=0$ if $s$ is degenerate and $\jac_{\mu^i}(s)=\frac{\pi}{|E|}$ if $s$ is a norm, where $|E|$ denotes the Lebesgue measure of John's ellipse of $\{v\in\R^2:s(v)\leq1\}$. Again, by identifying $(T_zM,g(z))$ with $(\R^2,|\cdot|)$ via a linear isometry, we are able to define the $\mu^i$-jacobian of a seminorm $s$ on $T_zM$.
\bd
    The \emph{inscribed Riemannian area} of $u\in N^{1,2}( U,X)$ is defined as
    $$\Area_{\mu^i}(u):=\int_U\jac_{\mu^i}(\apmd u_z)\;d\hm_g^2(z).$$
\ed

From \cite[Section~2.4]{LW16-energy-area} it follows that the inscribed Riemannian area and the parame-trized Hausdorff area are comparable; more explicitly 
\begin{align}\label{ineq:inscribed-Riem-and-Hausdorff-area}
    \frac{\pi}{4}\Area_{\mu^i}(u)\leq\Area(u)\leq\Area_{\mu^i}(u).
\end{align}

\bd
A map $u\in N^{1,2}( U,X)$ is called \emph{infinitesimally isotropic with respect to a Riemannian metric $g$ on $M$} if for almost every $z\in  U$ the approximate metric derivative $\apmd u_z$ is either zero or it is a norm and the John's ellipse of $\apmd u_z$ is a round ball with respect to $g$.
\ed

It follows from \cite{LW16-energy-area} that $\Area_{\mu^i}(u)\leq E^2_+(u,g)$ for all Riemannian metrics $g$ on $M$, with equality if and only if $u$ is infinitesimally isotropic. Moreover, by \cite{LW16-energy-area}, if $u\in N^{1,2}(U,X)$ is infinitesimally isotropic with respect to $g$, then it is \emph{infinitesimally $K$-quasiconformal with respect to $g$} with $K=\frac{4}{\pi}$ in the sense that
$$ (\rho_u(z))^2 \leq K\cdot \jac(\apmd u_z) $$
for almost every $z\in U$. 

Assume that $X$ is a complete metric space and $u\colon M\to X$ continuous and monotone. By arguing as in the proof of \cite[Proposition~3.5]{LW20-param}, we obtain that if $u\in N^{1,2}(M, X)$ and $u$ is infinitesimally $K$-quasiconformal with respect to $g$ then 
\begin{align}\label{eq:inf-qc-props}
    \MOD(\Gamma)\leq K\cdot\MOD(u\circ\Gamma)
\end{align}
for every family $\Gamma$ of curves in $M$. Conversely, if $u$ is a homeomorphism onto its image and satisfies \eqref{eq:inf-qc-props} then $u$ belongs to $N^{1,2}(M, X)$ and is infinitesimally $K$-quasiconformal (see \cite[Theorem 1.1]{Wil12}).\\

We define the \emph{trace} of a Sobolev map $u\in N^{1,2}( U,X)$ in the following way. Assume $ U\subset M\setminus\partial M$ is a Lipschitz domain. Then for every point $z\in\partial U$ there exists an open neighborhood $V\subset M$ and a biLipschitz mapping $\psi\colon(0,1)\times[0,1)\to M$ such that $\psi((0,1)\times(0, 1))=U\cap  V$ and $\psi((0,1)\times\{0\})=V\cap\partial U$. As $u$ is Sobolev, for almost every $s\in(0,1)$ the map $t\mapsto u\circ\psi(s,t)$ has an absolutely continuous representative which we denote by the same expression. The trace of $u$ 
$$tr(u)(\psi(s,0)):=\lim_{t\searrow 0}(u\circ\psi)(s,t)$$
is defined for almost every $s\in(0,1)$. It can be shown (see \cite[Section~1.12]{KS93}) that the trace is independent of the choice of the map $\psi$ and is in $L^2(\partial U,X)$. 

Assume that $M$ has $k\geq1$ boundary components and let $\Gamma$ be a disjoint union of $k$ Jordan curves in $X$. Recall that $[\Gamma]$ denotes the set of all weakly monotone parametrizations of $\Gamma$. Let $\Lambda(M,\Gamma,X)$ be the possibly empty family of Sobolev maps $u\in N^{1,2}(M,X)$ such that the trace $\trace(u)$ has a continuous representative in $[\Gamma]$.

\subsection{Relative $1$-Homotopy classes of Sobolev maps}\label{sec:prelimsRelHom}
We now introduce a notion of relative $1$-homotopy classes of Sobolev mappings; for more information we refer to \cite{SW20}. 

A finite collection $K$ of compact convex polytopes (called cells of $K$) in some $\R^n$ is a \emph{polyhedral complex} if each face of a cell is in $K$ and the intersection of two cells of $K$ is a face of each of them. We always equip $K$ with the induced metric from $\R^n$, implying that a $2$-cell $\Delta$ is isometric to a compact convex polygon in $\R^2$.

In the following let $M$ be a smooth surface with possibly non-empty boundary. A \emph{triangulation} of $M$ consists of a polyhedral complex $K$ and a homeomorphism $h\colon K\to M$, where $h$ restricted to any 2-cell $\Delta$ of $K$ is a $C^1$-diffeomorphism onto its image. The $j$-skeleton of $K$, denoted $K^j$, is the union of all cells of $K$ of dimension at most $j$. Let $\partial K\subset K^1$ be the preimage of $\partial M$ under $h$ and define $$\widehat{K}^1:=(K^1\setminus\partial K)\cup K^0.$$ 
Note that if $M$ has empty boundary, then $\widehat{K}^1=K^1$.

Consider a proper geodesic metric space $X$.
Two continuous maps ${\varrho,\varrho'\colon K^1\to X}$ are homotopic relative to a set $A\subset K^1$ if there exists a homotopy $H\colon K^1\times[0,1]\to X$ between $\varrho$ and $\varrho'$ satisfying $H(s,\cdot)=\varrho(s)=\varrho'(s)$ for every $s\in A$. If $A$ is empty, then relative homotopy agrees with the usual definition of homotopy. For $\varrho|_{\partial K},\varrho'|_{\partial K}\in[\Gamma]$ we say that $\varrho$ and $\varrho'$ are \emph{homotopic relative to $\Gamma$} in some ambient space $Y\supset X$, denoted $$\varrho\sim\varrho'\text{ rel }\Gamma\text{ in }Y,$$ if there exists a homotopy $H\colon K^1\times[0,1]\to Y$ between $\varrho$ and $\varrho'$ so that $H(\cdot,t)|_{\partial K}\in[\Gamma]$ for every $t\in[0,1]$. If $Y$ is not mentioned, we assume $X=Y$. The \emph{homotopy class of $\varrho$ relative to $\Gamma$} is the family
    $$[\varrho]_{\Gamma}:=\{\varrho'\colon K^1\to X: \varrho'\text{ continuous, }\varrho'|_{\partial K}\in[\Gamma],\varrho\sim\varrho'\text{ rel }\Gamma\}.$$
We will also use this notation if $\Gamma$ is empty; then the set $[\rho]_{\Gamma}$ coincides with the usual homotopy class $[\rho]$.

\bd
    An \emph{admissible deformation on a surface $M$} is a smooth map $\Phi\colon M\times\R^m\to M$, $m\in\N$, such that $\Phi_{\xi}:=\Phi(\cdot,\xi)$ is a diffeomorphism for every $\xi\in\R^m$ and $\Phi_0=\id_M$, and such that the derivative of $\Phi^p:=\Phi(p,\cdot)$ satisfies
    $$D\Phi^p(0)(\R^m)=\begin{cases} T_pM & \text{if } p\in\interior(M)\\
    T_p(\partial M) & \text{if } p\in\partial M.
    \end{cases}$$
\ed

A metric space $Y$ is called \emph{$\varepsilon$-thickening} of $X$, $\varepsilon>0$ if there exists an isometric embedding $\iota\colon X\to Y$ such that the Hausdorff distance between $\iota(X)$ and $Y$ is less than $\varepsilon$. For any compact metric space $X$ and any $\varepsilon>0$ there exists a $\varepsilon$-thickening $X_{\varepsilon}$ of $X$ that is again compact and satisfies a local quadratic isoperimetric inequality, i.e. there exist $C,l_0>0$ such that every Lipschitz curve $c\colon S^1\to X_{\varepsilon}$ of length $\ell(c)\leq l_0$ is the trace of a Sobolev map $u\in N^{1,2}(D, X_{\varepsilon})$ with
$$ \Area(u)\leq C\cdot\ell(c)^2.$$ This result follows from \cite{Wen08-sharp} and \cite[Lemma 3.3]{LWY20} if $X$ locally geodesic and from \cite[Lemma~5.1]{CF20} else. 
The definition of relative $1$-homotopy classes in \cite{SW20} highly depends on the existence of a quadratic isoperimetric inequality. Thus, for our definition of relative $1$-homotopy we will apply the statements of \cite{SW20} to the $\varepsilon$-thickening $X_{\varepsilon}$ of $X$ and use the following lemma to pass to a notion of relative 1-homotopy in $X$.

\bl\label{lem:Retraction}
    Let $X$ be a metric surface. Then there exists $\varepsilon>0$ and a continuous retraction $$R\colon X_{\varepsilon}\to X.$$
\el

In order to prove Lemma~\ref{lem:Retraction} we need some more notation. A metric space $X$ is an \emph{absolute neighbourhood retract (ANR)} if for each closed subset $A$ of a metric space $Y$, every continuous map $f\colon A\to X$ has a continuous extension $F\colon U\to X$ defined on some neighbourhood $U$ of $A$ in $Y$. By \cite[Corollary~14.8A]{Dav07}, every finite dimensional, locally contractible, compact metric space is an ANR. Thus, whenever $M$ is a smooth surface and $X$ a metric space homeomorphic to $M$, then both $M$ and $X$ are ANRs. 

\begin{proof}[Proof of Lemma \ref{lem:Retraction}]
    For every $n\in\N$ we consider $X$ as a subset of its $1/n$-thickening $X_{1/n}$ and define
    $$Y:=\bigsqcup_{n\in\N}X_{1/n}\Big/\sim,$$
    where for $x_i\in X_{1/n_i}$ we have that $x_1\sim x_2$ if and only if $d_Y(x_1,x_2)=0$ with
    $$d_Y(x_1,x_2):=\min_{y\in X}\left\{d_{X_{1/n_1}}(x_1,y)+d_{X_{1/n_2}}(x_2,y)\right\}.$$
    The $1/n$-thickening $X_{1/n}$ as well as $X$ embed isometrically into $(Y,d_Y)$ for every $n\in\N$. The inclusion map is in particular closed and by \cite[Proposition~14.2]{Dav07}, there is an open neighbourhood $U$ of $X$ in $Y$ and a continuous retraction $R\colon U\to X$. The lemma follows after choosing $\varepsilon:=1/n$ for $n\in\N$ large enough that $X_{\varepsilon}\subset U$ and restricting $R$ to $X_{\varepsilon}$. 
\end{proof}

Consider a metric surface $X$ and let $\varepsilon>0$ and $R\colon X_{\varepsilon}\to X$ be as in Lemma~\ref{lem:Retraction}. The $\varepsilon$-thickening $X_{\varepsilon}$ satisfies a quadratic isoperimetric inequality and we can thus apply the results form \cite{SW20}. Let $\Phi\colon M\times\R^m\to M$ an admissible deformation on $M$, whose existence follows from \cite[Proposition~3.2]{SW20}. For a triangulation $h\colon K\to M$ of $M$ and $\xi\in\R^m$ denote by $h_{\xi}\colon K\to M$ the triangulation given by $h_{\xi}:=\Phi_{\xi}\circ h$. Furthermore, for $\xi\in\R^m$ and $u\in N^{1,2}(M,X)$ we denote by $u\circ h_{\xi}|_{K^1}$ the map agreeing with $u\circ h_{\xi}$ on $K^1\setminus\partial K$ and with $\trace(u)\circ h_{\xi}$ on $\partial K$. Fix a Riemannian metric $g$ on $M$. In \cite[Section~3]{SW20} it is shown that for every $u\in\Lambda(M,\Gamma,X)$ and every triangulation $h\colon K\to M$ of $M$ there exists a ball $B_{\Phi,h}\subset\R^m$ centered at the origin such that for almost all $\xi,\zeta\in B_{\Phi,h}$ the maps $u\circ h_{\xi}|_{K^1}$ and $u\circ h_{\zeta}|_{K^1}$ have continuous representatives which are homotopic relative to $\Gamma$ in $X_{\varepsilon}$. After postcomposition with $R$, the continuous representatives of $u\circ h_{\xi}|_{K^1}$ and $u\circ h_{\zeta}|_{K^1}$ are homotopic relative to $\Gamma$ in $X$. We denote the common relative homotopy class by $u_{\#,1}[h]$. Note that $u_{\#,1}[h]$ is independent of the choice of deformation $\Phi$ and inducing the same relative homotopy class is independent of the triangulation $h$ (see \cite[Theorem~4.1]{SW20}). Moreover, if $u$ is continuous, then $u_{\#,1}[h]=[u\circ h|_{K^1}]_{\Gamma}$ for every triangulation $h$ of $M$.

Two maps $u,v\in\Lambda(M,\Gamma,X)$ are \emph{$1$-homotopic relative to $\Gamma$}, denoted $u\sim_1 v\rel\Gamma,$ if $$u_{\#,1}[h]=v_{\#,1}[h]$$ for one and thus any triangulation $h$ of $M$.

\subsection{Injective metric spaces}\label{Sec:InjectiveMetricSpaces}

A metric space $E$ is \emph{injective} if for every metric space $Z$, any $1$-Lipschitz map $A\to E$, defined on a subset $A\subset Z$, extends to a $1$-Lipschitz map $Z\to E$. For any metric space $X$ there exists an injective metric space $E(X)$, called \emph{injective hull of $X$}, which contains $X$ and which is minimal in an appropriate sense among injective metric spaces containing $X$, see \cite{Isb64}. Note that $E(X)$ is unique up to isometry and if $X$ is compact then so is $E(X)$. Moreover, if $X\subset Y$ then $E(X)$ can be considered as a subset of $E(Y)$. For $r>0$, we denote by $N_{r}(X)$ the $r$-neighbourhood of $X$ in $E(X)$. 

\bl\label{lem:retraction-injective-hull}
    Let $X$ be a metric surface. Then there exists $\varepsilon, r_0>0$ and a continuous retraction $$R_0\colon (N_{r_0}(X))_{\varepsilon}\to X,$$ 
    where $(N_{r_0}(X))_{\varepsilon}$ denotes the $\varepsilon$-thickening of $N_{r_0}(X)$.
\el

Lemma~\ref{lem:retraction-injective-hull} follows by arguing exactly as in the proof of Lemma~\ref{lem:Retraction} and will be used to apply results from \cite{SW20} to mappings with image in neighborhood $N_{r_0}(X)$.

\section{Non-trivial Sobolev maps}
\label{sec:existence}

The purpose of this section is to establish the existence of a Sobolev map in $\Lambda(M,\varphi,X)$. The following variant of \cite[Proposition~5.1]{MW21} will play a crucial role.
\bp\label{prop:Lipschitz-map-MW21}
    Let $X$ be a locally geodesic metric space homeomorphic to a smooth surface $M$. Suppose $\Omega\subset M$ is biLipschitz equivalent to $\overline{D}$ and $J\subset X$ homeomorphic to $\overline{D}$. Consider a biLipsichitz map $\chi\colon I\to\partial J$ for $I\subset\partial\Omega$ connected. If $\hm^2(J)<\infty$ and $\ell(\partial J)<\infty$ then there exists a constant $c>0$ with the following property. For every $r>0$ there is a Lipschitz map $v\colon \Omega\to E(J)$ with $\Area(v)\leq c$ and such that $v|_{\partial\Omega}$ parametrizes $\partial J$, $v|_I=\chi$ and $\Image(v)\subset N_{r}(J)$.
\ep
\begin{proof}
    In the special case that $J$ equipped with the subspace metric is geodesic, the statement follows from \cite[Proposition~5.1]{MW21}. Note that the parametrization of $v|_I$ can be prescribed by gluing a Lipschitz homotopy of zero area along $I$ and reparametrizing. For the existence of such a Lipschitz homotopy consider e.g. \cite[Proposition~3.6]{LWY20}. 
    
    For the general case, let $d$ be the metric on $X$ and denote by $d_{J}$ the length metric on $J$. The identity map $\pi\colon Y:=(J,d_{J}) \to (J, d)$ is a $1$-Lipschitz homeomorphism which preserves lengths of curves and the Hausdorff $2$–measure of Borel subsets, compare to \cite[Lemma~2.1]{LW20-param}. In particular, the metric space $Y$ is a geodesic space homeomorphic to $\overline{D}$ with rectifiable boundary and finite Hausdorff $2$-measure. Moreover, as $\chi(I)\subset\partial J$ is a biLipschitz curve, the inverse of $\pi$ restricted to $\chi(I)$ is Lipschitz. The special case implies the existence of a constant $c>0$ with the following property. For every $r>0$ there is a Lipschitz map $v\colon\Omega\to E(Y)$ with $\Area(v)\leq c$ and such that $v|_{\partial\Omega}$ parametrizes $\partial Y$, $v|_{I}=\pi^{-1}\circ\chi|_{I}$ and $\Image(v)\subset N_{r}(Y)$. By injectivity of $E(J)$, the $1$-Lipschitz map $\pi$ extends to a $1$-Lipschitz map $\overline{\pi}\colon E(Y)\to E(J)$. The Lipschitz map $\overline{v}:=\overline{\pi}\circ v$ fulfills the desired properties. 
\end{proof}
 
In order to apply Proposition~\ref{prop:Lipschitz-map-MW21}, we decompose $M$ and $X$ into Jordan domains, using a similar strategy as in \cite{FM21}. A surface of genus $0$ is called \emph{cylinder} if it has two boundary components and \emph{Y-piece} if it has three boundary components.

\bl\label{lem:decomp}
    Let $X$ be a locally geodesic metric space homeomorphic to a smooth surface $M$ and let $\varphi\colon M\to X$ be a homeomorphism. Then there exists a triangulation $h\colon K\to M$ of $M$ and a biLipschitz embedding $\chi\colon h(\widehat{K}^1)\to X$ such that $\varphi\circ h|_{\widehat{K}^1}$ and $\chi\circ h|_{\widehat{K}^1}$ are homotopic relative to $K^0\cap\partial K$ and $h(\Delta)$ is biLipschitz equivalent to $\overline{D}$ and $h(\Delta)\cap\partial M$ is connected for any 2-cell $\Delta$ of $K$.
\el

Recall that $\widehat{K}^1:=(K^1\setminus\partial K)\cup K^0$, which agrees with $K^1$ if $\partial K$ is empty.

\begin{proof}
    If $M$ is a disc, there's nothing to show. If $M=S^2$, equipped with the standard metric on $S^2$, we can choose three simple closed geodesics decomposing $M$ into eight domains $\Omega_{k,l}\subset M$, $k\in\{1,2,3,4\}$, $l\in\{1,2\}$, that are triangles on the sphere and biLipschitz equivalent to $\overline{D}$. If $M$ is not of disc- or sphere-type, depending on its topology, endow $M$ with a hyperbolic or flat Riemannian metric. Choose a collection of simple closed geodesics decomposing $M$ into smooth Y-pieces and cylinders $M_k$ that intersect at most one boundary component of $M$. It is a standard result from hyperbolic geometry that if $M_k$ is a Y-piece, then it is isometric to the partial gluing of the boundary of two right-angled hexagons $\Omega_{k,1}, \Omega_{k,2}\subset\mathbb{H}$, see e.g. \cite[Proposition 3.1.5]{Bus10}. Whenever $M_k$ is a cylinder, then a similar decomposition into isometric rectangles $\Omega_{k,1}, \Omega_{k,2}\subset\R^2$ is possible. In either case $\Omega_{k,1}$ and $\Omega_{k,2}$ are biLipschitz equivalent to the closed unit disc $\overline{D}$. 
    
    Fix a triangulation $h\colon K\to M$ of $M$ such that every $\Omega_{k,l}$ is the image of a 2-cell $\Delta_{k,l}$ of $K$ under $h$. In particular, it has to hold that $h|_{\Delta_{k,l}}\colon\Delta_{k,l}\to\Omega_{k,l}$ is a $C^1$-diffeomorphism, which is only possible if $\Delta_{k,l}$ is the same type of polytope as $\Omega_{k,l}$. For every edge $e_i\in\widehat{K}^1$ we set $\alpha_i:=h(e_i)$ and denote by $a_i^1,a_i^2$ the endpoints of $\alpha_i$. There are piecewise geodesic biLipschitz curves $\beta_i$ in $X$ arbitrary close to $\varphi(\alpha_i)$ and with endpoints $\varphi(a_i^1),\varphi(a_i^2)$, see \cite[Lemma~4.2]{LW20-param}. By arguing as in the proof of \cite[Lemma~3.1]{FM21}, we can modify each $\beta_i$ in an arbitrary small neighbourhood of $\beta_i$, while fixing an endpoint $\varphi(a_i^j)$ if $a_i^j\in\partial M$, such that $\bigcup_i \beta_i$ is biLipschitz equivalent to $\bigcup_i\alpha_i=h(\widehat{K}^1)$. In particular, there exists a biLipschitz map $\chi\colon\bigcup_i\alpha_i\to\bigcup_i\beta_i$ sending each $\alpha_i$ to $\beta_i$ and a homotopy relative to $K^0\cap\partial K$ between $\varphi\circ h|_{\widehat{K}^1}$ and $\chi\circ h|_{\widehat{K}^1}$.
\end{proof}

With Proposition~\ref{prop:Lipschitz-map-of-bdd-area-in-Xeps} and Lemma~\ref{lem:decomp} at hand, we are able to prove the following generalization of \cite[Proposition~5.1]{MW21}.

\bp\label{prop:Lipschitz-map-of-bdd-area-in-Xeps}
    Let $X$ be a locally geodesic metric space homeomorphic to a smooth surface $M$ and $\varphi\colon M\to X$ a homeomorphism.
    If $\hm^2(X)<\infty$ and $\ell(\partial X)<\infty$ then there exists a constant $C>0$ with the following property. For every $r\in(0,r_0]$ there exists a Lipschitz map $v\colon M\to E(X)$ with $\Area(v)\leq C$ and such that
    \begin{enumerate}
        \item $v|_{\partial M}$ parametrizes $\partial X$,
        \item $v|_{h(\widehat{K}^1)}=\chi|_{h(\widehat{K}^1)}$,
        \item $\Image(v)\subset N_{r}(X)$,
        \item $v\sim_1\varphi\rel\partial X$ in $(N_{r_0}(X))_{\varepsilon}$,
    \end{enumerate}
    where $h\colon K\to M$, $\chi\colon M\to X$ are as in Lemma~\ref{lem:decomp} and $\varepsilon, r_0>0$ are as in Lemma~\ref{lem:retraction-injective-hull}.
\ep

\begin{proof}
    Denote by $\Delta_1,...,\Delta_N$ the 2-cells of $K$ and define $\Omega_i:=h(\Delta_i)$ and $I_i:=h(\widehat{K}^1)\cap\partial\Omega_i$. Let $J_i$ be the closure of the Jordan domain obtained by cutting along $\chi(I_i)$. By Proposition~\ref{prop:Lipschitz-map-MW21}, for each $i$ there exists a constant $c_i>0$ such that the following holds. For every $r>0$ there is a Lipschitz map $v_i\colon \Omega_i\to E(J_i)$ with $\Area(v_i)\leq c_i$ and such that $v_i|_{\partial\Omega}$ parametrizes $\partial J_i$, $v_i|_{I_i}=\chi|_{I_i}$ and $\Image(v_i)\subset N_{r}(J)$. After gluing all $v_i$ along corresponding boundaries we obtain a Lipschitz map $v\colon M\to E(X)$ of area less than $C:=\sum_{i=1}^Nc_i$ satisfying (i), (ii) and (iii), where the constant $C$ is independent of $r$.
    
    By Lemma \ref{lem:decomp}, the maps $\varphi\circ h|_{\widehat{K}^1}$ and $\chi\circ h|_{\widehat{K}^1}=v\circ h|_{\widehat{K}^1}$ are homotopic relative to $K^0\cap\partial K$. Let $e\subset\partial K$ be a $1$-cell of $K$. As $\varphi\circ h|_e$ and $v\circ h|_e$ both parametrize $\varphi(h(e))\subset\partial X$, there exists a homotopy $F_e\colon e\times[0,1]\to X$ relative endpoints between $\varphi\circ h|_e$ and $v\circ h|_e$ such that $F_e(\cdot,t)$ parametrizes $\varphi(h(e))$ for every $t\in[0,1]$. After gluing the obtained homotopies along corresponding points in $K^0\cap\partial X$, we obtain a homotopy $H\colon K^1\times[0,1]\to X$ between $v\circ h|_{K^1}$ and $\varphi\circ h|_{K^1}$ such that $H(\cdot,t)|_{\partial K}$ is a parametrization of $\partial X$ for every $t\in[0,1]$.
\end{proof}

The proof of Theorem~\ref{thm:existence-Sobolev} uses similar arguments as in \cite[Section~8]{FW21} and the proofs of \cite[Proposition~6.1]{SW20} and \cite[Theorem~1.4]{MW21}. Let $X$ be a complete metric space and $M$ a smooth surface. A family $\mathcal{F}$ of continuous maps from $M$ to $X$ is said to \emph{satisfy the condition of cohesion} if there exists $\eta>0$ such that each $u\in\mathcal{F}$ satisfies $\ell(u\circ c)\geq\eta$ for every non-contractible closed curve $c$ in $M$. Furthermore, two $1$-cells in $K$ are called \emph{non-neighbouring} whenever they do not intersect.

\begin{proof}[Proof of Theorem \ref{thm:existence-Sobolev}]
Let $\varepsilon, r_0>0$ and $R_0\colon (N_{r_0}(X))_{\varepsilon}\to X$ be as in Lemma~\ref{lem:retraction-injective-hull}.
Consider the triangulation $h\colon K\to M$ and the biLipschitz map $\chi\colon h(\widehat{K}^1)\to X$ from Lemma~\ref{lem:decomp}. By Proposition~\ref{prop:Lipschitz-map-of-bdd-area-in-Xeps}, there exists a constant $C>0$ and a sequence of Lipschitz mappings $(v_n)$ such that each $v_n\colon M\to E(X)$ has area less than $C$ and satisfies properties (i) to (iv) from  Proposition~\ref{prop:Lipschitz-map-of-bdd-area-in-Xeps} with $r:=1/n$.

As $\chi$ is biLipschitz, it holds that the constant
$$ \eta:=\inf\{\dist(\chi(h(e_1)),\chi(h(e_2))):e_1,e_2\text{ are non-neighbouring 1-cells of } K \}$$
is strictly positive. Observe that if $c$ is a non-contractible closed curve in $M$, then $h^{-1}\circ c$ intersects at least two non-neighbouring 1-cells $e_1,e_2$ of $K$. Denote the intersection point of $c$ and $e_i$ by $p_i$. As $X$ embeds isometrically into $E(X)$ it holds that $$\ell(v_n\circ c)\geq d_{E(X)}(v_n(h(p_1)),v_n(h(p_2q)))=d(\chi(h(p_1)),\chi(h(p_2)))\geq\eta.$$
Thus, $(v_n)$ satisfies the condition of cohesion.
By Morrey's $\varepsilon$-conformality lemma (see \cite[Theorem~1.2]{FW20}) and inequality (\ref{ineq:inscribed-Riem-and-Hausdorff-area}), we obtain a sequence of hyperbolic (or flat) metrics $(g_n)$ on $M$ such that
$$E_+^2(v_n,g_n)\leq\frac{4}{\pi}\Area(v_n)+1\leq\frac{4}{\pi}C+1.$$
By \cite[Proposition~8.4]{FW21}, there exists a uniform constant $\varepsilon>0$ such that the relative systole of $(M,g_n)$ is bounded from below by $\varepsilon$. Hence, we can apply the Mumford compactness theorem to obtain orientation preserving diffeomorphisms $\phi_n\colon M\to M$ such that a subsequence of $(\phi_n^*g_n)$ converges to a hyperbolic (or flat) metric $h$ on $M$ (see \cite[Theorem~3.3]{FW21} and e.g. \cite[Theorem~4.4.1]{DHT10} for the fact that the diffeomorphisms may be chosen to be orientation preserving). This convergence implies that for maps $w_n:=v_n\circ\phi_n\in\Lambda(M,\partial X,E(X))$ it holds that
$$E_+^2(w_n,h)\leq C_n\cdot E_+^2(v_n,g_n),$$
where $C_n\geq1$ tends to $1$ as $n\to\infty$. By the Rellich-Kondrachov compactness theorem (see \cite[Theorem~1.13]{KS93}), a further subsequence of $(w_n)$ converges in $L^2(M,E(X))$ to a finite energy map $w$. As each $v_n$ is $1$-homotopic to $\varphi$ relative to $\partial X$ in $(N_{r_0}(X))_{\varepsilon}$ and every $\phi_n$ is orientation preserving, the maps $w_n$ induce the same orientation on $\partial X$. From \cite[Theorem~4.7]{SW20} we know that there exists $j_0\in\N$ such that for every $j\geq j_0$ the map $w_{n_j}$ is $1$-homotopic to $w_{n_{j_0}}$ relative to $\partial X$ in $(N_{r_0}(X))_{\varepsilon}$. Then, the maps $u_j:= w_{n_j}\circ\phi_{n_{j_0}}^{-1}\in\Lambda(M,\partial X,E(X))$ satisfy for $j\geq j_0$
$$u_j\sim_1w_{n_{j_0}}\circ\phi_{n_{j_0}}^{-1}=v_{n_{j_0}}\sim_1\varphi\text{ rel }\partial X\text{ in }(N_{r_0}(X))_{\varepsilon}. $$
The sequence $(u_j)$ converges in $L^2(M,X)$ to the map $u:= w\circ\phi_{n_{j_0}}^{-1}$ and we set $g:=(\phi_{n_{j_0}}^{-1})^*h$. Since the image of $u_j$ is in $N_{1/n_j}(X)$ we have that the essential image of $u$ is in $X$ and thus $u$ can be viewed as an element of $N^{1,2}(M,X)$. This already shows that $u\in\Lambda(M,\partial X,X)$ if $\partial X=\emptyset$. Assume now that $\partial X$ is not empty. By \cite[Proposition~8.3]{FW21}, the sequence $(\trace(u_j))$ is equicontinuous as $(u_j)$ still satisfies the condition of cohesion. Hence, by Arzelà-Ascoli, a subsequence of $(\trace(u_j))$ converges uniformly to some continuous map $\gamma\colon\partial M\to X$, which is a weakly monotone parametrization of $\partial X$. The convergence of $(\trace(u_j))$ in $L^2(\partial M,X)$ to $\trace(u)$, see \cite[Theorem~1.12.2]{KS93}, implies that $\trace(u)$ coincides with $\gamma$. This shows that $u\in\Lambda(M,\partial X,X)$. It follows from \cite[Theorem~4.7]{SW20} that the map $u$ is $1$-homotopic to $\varphi$ relative to $\partial X$ in $(N_{r_0}(X))_{\varepsilon}$. Since both $u$ and $\varphi$ have image in $X$ and any homotopy in $(N_{r_0}(X))_{\varepsilon}$ can be projected by $R_0$ to a homotopy with image in $X$, $u$ is in fact $1$-homotopic to $\varphi$ relative to $\partial X$ in $X$.
\end{proof}

Note that we exclude the case of a sphere in Theorem~\ref{thm:existence-Sobolev}, since for spaces with non-contractible universal coverings we can not apply the same methods to find a converging subsequence of $(v_n)$. The methods used to prove the existence of energy minimizing spheres in the smooth case (see \cite{SU81}) have not yet been extended to this generality.

In a next step we apply a direct variational method to show the existence of an energy minimizing pair in $\Lambda_{\text{metr}}(M,\varphi,X)$.

\bt\label{thm:Sobolev-energy-min}
    Let $M$ be a smooth surface that is not a sphere, $X$ a metric space homeomorphic to $M$ and $\varphi\colon M\to X$ a homeomorphism. If $\Lambda(M,\varphi,X)$ is not empty, then there exists an energy minimizing pair $(u,g)\in\Lambda_{\text{metr}}(M,\varphi,X)$.    
\et

\begin{proof}
    Take an energy minimizing sequence $(v_n,g_n)$ in $\Lambda_{\text{metr}}(M,\varphi,X)$, i.e. a sequence of pairs $(v_n,g_n)\in\Lambda_{\text{metr}}(M,\varphi,X)$ satisfying 
    $$E_+^2(v_n,g_n)\to\inf\{E_+^2(v,g):(v,g)\in\Lambda_{\text{metr}}(M,\varphi,X)\}$$
    as $n$ tends to infinity. Every non-contractible closed curve $c$ in $X$ satisfies $\ell(c)\geq\eta$ for some $\eta>0$ as $X$ is homeomorphic to a smooth surface. Thus, by arguing as in the proofs of Propositions~8.4 and 8.3 in \cite{FW21}, we obtain that the relative systole of $(M,g_n)$ is bounded away from zero independently of $n$ and the sequence $(\trace(v_n))$ is equicontinuous. Note that the maps in \cite[Section 8]{FW21} were additionally assumed to be continuous, but the proofs of Propositions~8.4 and 8.3 in \cite{FW21} can be adapted to the current setting. We proceed as in the proof of Theorem~\ref{thm:existence-Sobolev} to obtain a map $u\in\Lambda(M,\varphi,X)$ and a hyperbolic metric $g$ on $M$ with the following property. The map $u$ is $1$-homotopic to $\varphi$ relative to $\partial X$ and after precomposing each $v_n$ with a suitable diffeomorphism of $M$ and passing to a subsequence, the maps $v_n$ converge to $u$ in $L^2(M,X)$. The statement follows from lower semicontinuity of energy.
\end{proof}

\section{Continuity of energy minimizers}\label{sec:continuity}

The goal of this section is to provide a proof of Theorem~\ref{thm:continuity-of-energy-min-Sobolev}, equipping us with the right regularity of energy minimizers. For an arbitrary map $v\colon M\to X$ from a smooth surface $M$ to a metric space $X$ and for $z\in M$ and $\delta>0$ the \emph{essential oscillation of $v$ in the $\delta$-ball around $z$} is defined by 
$$\osc(v,z,\delta):= \inf\{\diam(v(A)): \text{$A\subset M\cap B(z,\delta)$ subset of full measure}\}.$$

We can show as in the proof of \cite[Theorem~1.3]{MW21} that Theorem~\ref{thm:continuity-of-energy-min-Sobolev} is implied by the following generalization of \cite[Proposition~4.1]{MW21}.

\bp\label{prop:essential-oscillation-energy-min}
    Let $X$ be a locally geodesic metric space homeomorphic to a smooth surface $M$ with non-empty boundary and $\varphi\colon M\to X$ a homeomorphism. If $(u,g)\in\Lambda_{\text{metr}}(M,\varphi,X)$ is an energy minimizing pair, then for every $\varepsilon>0$ there exists $\delta>0$ such that $\osc(u,z,\delta)<\varepsilon$ for every $z\in M$.
\ep

The proof of Proposition~\ref{prop:essential-oscillation-energy-min} is very similar to the proof of \cite[Proposition~4.1]{MW21}. We repeat the most important steps and adapt them to the current setting. For details we refer to \cite[Section~4]{MW21}.

Consider $M$, $X$ and $\varphi\colon M\to X$ as in Proposition~\ref{prop:essential-oscillation-energy-min}. Let $(u,g)\in\Lambda_{\text{metr}}(M,\varphi,X)$ be an energy minimizing pair and let $\varepsilon>0$. For an application of the Courant-Lebesgue Lemma, we consider a family $\mathcal{F}$ of $2$-biLipschitz mappings from $\overline{D}$ to $M$ such that every point in $M$ is contained in the image of at least one map in $\mathcal{F}$; compare to the usage of Courant-Lebesgue in the proof of \cite[Proposition~4.1]{FM21}. Fix $z\in M$ and choose $\psi\in\mathcal{F}$ with $z\in\Image(\psi)$. Let $\delta>0$ be the constant from the Courant-Lebesgue Lemma applied to the map $u\circ\psi$. Note that we can choose $\delta>0$ so small that for any $r\in(\delta,\sqrt{\delta})$ the set $\psi(\partial B(\psi^{-1}(z),r)\cap\overline{D})$ is contractible in $M$ and intersects at most one boundary component of $M$. As in the proof of \cite[Proposition~4.1]{MW21}, we find
$r\in(\delta,\sqrt{\delta})$ such that for $W:=\psi(D\cap B(\psi^{-1}(z),r))$ we have that the image of the trace of $u|_W$ is contained in a Jordan domain $\Omega$ of diameter less than $\varepsilon$ that is bounded by a bi Lipschitz Jordan curve or the concatenation of a biLipschitz Jordan arc and a connected subcurve of $\partial X$, compare to \cite[Lemma~4.2]{MW21}.

We are left to argue why the set $N:=\{w\in W:u(w)\in X\setminus\overline{\Omega}\}$ is negligible (compare to \cite[Lemma~4.3]{MW21}). For this we need the following lemma.

\bl \label{lem:Lip-retract}
    Let $X$ be a locally geodesic metric space homeomorphic to a smooth surface $M$ with non-empty boundary. If $\Omega\subset X$ is a Jordan domain that is either bounded by a biLipschitz Jordan curve or by the concatenation of a connected subcurve of $\partial X$ and a biLipschitz Jordan arc, then there exists a Lipschitz retraction $\varrho\colon X\to \overline{\Omega}$.
\el

\begin{proof}
    We provide a proof for $\partial\Omega$ being a Jordan curve contained in the interior of $X$, the case of $\partial\Omega$ intersecting $\partial X$ only needing minor adaptions in the following arguments. To get a better understanding of the construction described below, consider Figure~\ref{fig:LipschitzRetraction}.

    \begin{figure}
        \centering
        \includegraphics[width=.6\textwidth]{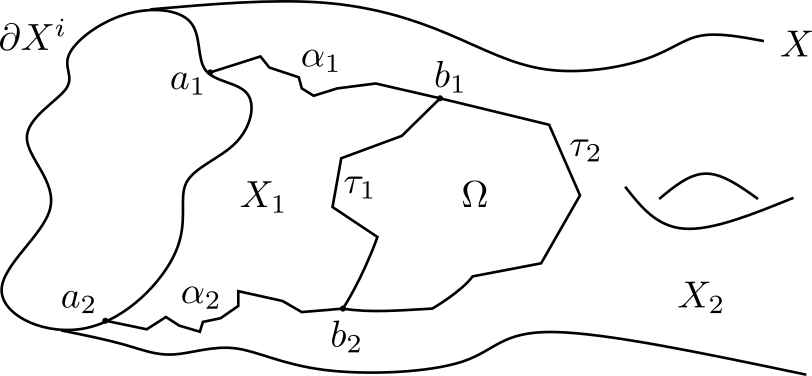}
        \caption{Construction of Lipschitz Retraction in proof of Lemma \ref{lem:Lip-retract}.}
        \label{fig:LipschitzRetraction}
    \end{figure}
    
    Choose a boundary component $\partial X^i$ and two points $a_1,a_2\in\partial X^i$. Decompose $\partial\Omega$ into two subcurves $\tau_1,\tau_2$ with endpoints $\tau_j(0)=:b_1$ and $\tau_j(1)=:b_2$. From \cite[Lemma 4.2]{LW20-param} we obtain the existence of two disjoint biLipschitz curves $\alpha_k$ connecting $a_k$ and $b_k$ in such a way that the concatenation of $\alpha_1$, $\alpha_2$, $\tau_1$ and  the corresponding component of $\partial X^i\setminus\{a_1,a_2\}$ bounds a Jordan domain $X_1$. Similarly, define $X_2$ to be the set bounded by the concatenation of $\alpha_1$, $\alpha_2$, $\tau_2$ and  the other component of $\partial X^i\setminus\{a_1,a_2\}$. After possibly changing $\alpha_k$ in the vicinity of $b_k$ as in the proof of \cite[Lemma~3.1]{FM21} and redefining $b_k$ and $\tau_j$, we can assume that the concatenation of $\alpha_1$, $\alpha_2$ and $\tau_j$ is biLipschitz.
    
    Let $\varrho_j$ be the Lipschitz map agreeing with the identity on the image of $\tau_j$ and sending every point on the curve $\alpha_k$ to $b_k$. Since $\tau_j$ is a biLipschitz curve, we can apply McShane's theorem to obtain the existence of a Lipschitz map $\overline{\varrho}_j\colon X_j\to \tau_j$ extending $\varrho_j$. Then, the map $\varrho\colon X\to \overline{\Omega}$ agreeing with $\overline{\varrho}_j$ on $X_j$ and with the identity on $\Omega$ is Lipschitz, as the intersections of respective domains are biLipschitz curves.
\end{proof}

Recall that $(u,g)\in\Lambda_{\text{metr}}(M,\varphi,X)$ is an energy minimizing pair. By \cite[Corollary~1.3]{FW20}, the map $u$ is infinitesimally isotropic with respect to $g$ and, by \cite[Proposition~1.1]{FW20}, $u$ minimizes the inscribed Riemannian area $\Area_{\mu^i}(u)$ among all maps in $\Lambda(M,\varphi, X)$.

We want to show that $N$ is negligible and suppose to the contrary that $N$ is not negligible. As in the proof of \cite[Lemma~4.3]{MW21}, a Fubini-type argument implies that $\Area_{\mu^i}(u|_N)>0$. Let $\varrho\colon X\to \overline{\Omega}$ be a Lipschitz retract as in Lemma~\ref{lem:Lip-retract} and denote by $v$ the map agreeing with $u$ on $M\setminus W$ and with $\varrho\circ u$ on $W$. Since the image of the trace of $u|_W$ is contained in $\Omega$, it follows from the Sobolev gluing theorem \cite[Theorem 1.12.3]{KS93} that $v\in N^{1,2}(M, X)$ and $\trace(v) = \trace(u)$. Moreover, it holds that $v\in\Lambda(M,\varphi,X)$. Indeed, if $h\colon K\to M$ is a triangulation of $M$ with $W\subset h(\text{int}(\Delta))$ for some 2-cell $\Delta$ of $K$ and $\Phi$ is any admissible deformation on $M$. Then for every sufficiently small $\xi$ the image of $K^1$ under $h_{\xi}$ does not intersect $W$ and therefore $u\circ h_{\xi}|_{K^1}=v\circ h_{\xi}|_{K^1}$.

Observe that $\Area_{\mu^i}(v|_N)=0$ as $v(N)\subset \partial\Omega$ and hence, $$\Area_{\mu^i}(u)=  \Area_{\mu^i}(u|_{M\setminus N}) + \Area_{\mu^i}(u|_N)> \Area_{\mu^i}(u|_{M\setminus N}) = \Area_{\mu^i}(v),$$
contradicting the area minimization property of $u$. Hence, $N$ is negligible and we finished the sketch of proof of Proposition~\ref{prop:essential-oscillation-energy-min}.

\section{Almost homeomorphism}\label{sec:monotone}

In this section we provide a proof of Theorem~\ref{thm:cont-and-inf-iso-implies-monotone} for surfaces that are not homeomorphic to a disc. If $M$ is a disc, we refer to the proof of \cite[Theorem~1.2]{LW20-param}.

A compact metric space is called \emph{cell-like} if it admits an embedding into the Hilbert cube in which it is null-homotopic in every neighborhood of itself. A continuous surjection $v\colon Y\to Z$ between metric spaces is called \emph{cell-like} if $v^{-1}(z)$ is cell-like, and in particular compact, for every $z\in Z$. Cell-like mappings are closely related to uniform limits of homeomorphisms as illustrated by the following theorem of Moore (see e.g. \cite[p.~116]{Edw78} or \cite[Theorem~25.1]{Dav07} for closed surfaces and \cite[Theorem~A]{Sie72} for compact surfaces with non-empty boundary).
\bt \label{thm:Moore}
    Let $M$ be a smooth surface and $v\colon M\to X$ a cell-like map such that $v|_{\partial M}\colon\partial M\to\partial X$ is cell-like if $\partial M$ is non-empty. Then $X$ is homeomorphic to $M$ and $v$ is a uniform limit of homeomorphisms.
\et
Moreover, by arguing exactly as in the proof of \cite[Proposition~2.9]{LW20-param} and using the fact that every metric surface is an ANR, we obtain the next proposition.
\bp
\label{prop:genOf2.9}
Let $M$ be a smooth surface with possibly non-empty boundary and let $X$ be a metric surface homoeomorphic to $M$. If $v\colon M\to X$ is a continuous surjection, then the following statements are equivalent:
\begin{enumerate}
    \item $v$ is monotone,
    \item $v$ is cell-like,
    \item $v$ is a uniform limit of homeomorphisms $v_i\colon M\to X$.
\end{enumerate}
\ep

From now on, we assume that $X$ is a locally geodesic metric space homeomorphic to a smooth surface $M$ with non-empty boundary.
The next theorem generalizes \cite[Theorem~4.1]{LW20-param} and will play a crucial role in the proof of Theorem~\ref{thm:cont-and-inf-iso-implies-monotone}.

\bt
\label{thm:cell-like}
Let $v\colon M\to X$ be a continuous surjection satisfying the following properties:
\begin{enumerate}
    \item The restriction of $v$ to $\partial M$ is a weakly monotone parameterization of $\partial X$. 
    \item Whenever $T\subset X$ is a single point or biLipschitz homeomorphic to a closed interval, every connected component of $v^{-1}(T)$ is cell-like.
\end{enumerate}
Then $v$ is a cell-like map.
\et

\begin{proof}
    We follow the same strategy as in the proof of \cite[Theorem~4.1]{LW20-param}.
    By the monotone-light factorization theorem due to Eilenberg and Whyburn (see e.g. \cite[Theorem~3.5]{You51}) there exists a compact metric space $Z$ and continuous surjective maps $v_1\colon M\to Z$, $v_2\colon Z\to X$, where $v_1$ is monotone and $v_2$ is light, such that for every $z\in Z$ the fibers $v_1^{-1}(z)$ are the connected components of $v^{-1}(v_2(z))$.  Recall that the map $v_2\colon Z\to X$ is light if $v_2^{-1}(x)$ is totally disconnected for every $x\in X$. The mappings $v_1$ and $v_1|_{\partial M}$ are cell-like as $v$ satisfies (i) and (ii). Hence, it follows from Theorem~\cite[Theorem~A]{Sie72} that $Z$ is homeomorphic to $M$ and $v_1$ is a uniform limit of homeomorphisms. Identify $Z$ with $M$. Since \cite[Theorem~1.4]{Lac69} holds for an arbitrary ANR, we can follow as in the proof of \cite[Lemma~4.4]{LW20-param} that the map $v_2$ satisfies properties (i) and (ii). If $v_2$ is cell-like, then so is $v$. Thus, it suffices to consider the case where $v$ is in addition a light map.
    
    As an application of Theorem~\ref{thm:Moore}, Lemma~2.8 in \cite{LW20-param} remains true for an arbitrary surface $(M,g)$ instead of $(S^2,g_{\mathrm{Eucl}})$. Therefore, Lemma~4.5 and Lemma~4.6 in \cite{LW20-param} can be generalized to our setting. Proceeding as in the last paragraph of \cite[Section~4]{LW20-param} completes the proof of Theorem~\ref{thm:cell-like}.
\end{proof}

The rest of this section is devoted to the proof of Theorem~\ref{thm:cont-and-inf-iso-implies-monotone}. Recall that $M$ has non-empty boundary and we excluded the case of $M$ being a disc. Let $\varphi\colon M\to X$ be a homeomorphism. Consider an energy minimizing pair $(u,g)\in\Lambda_{\text{metr}}(M,\varphi,X)$. As above, the map $u$ minimizes the inscribed Riemannian area $\Area_{\mu^i}$ among all maps in $\Lambda(M,\varphi,X)$ with respect to $g$.
By Proposition~\ref{prop:genOf2.9}, it suffices to show that the hypotheses of Theorem~\ref{thm:cell-like} hold. The map $u$ is surjective for topological reasons and $u\in\Lambda(M,\partial X, X)$ implies that $u$ satisfies (i). Towards a contradiction, assume there exists $T$ as in (ii) such that some connected component $K$ of $u^{-1}(T)$ is not cell-like. Let $\varepsilon>0$ be so small that the $\varepsilon$-neighbourhood of $T$ in $X$ is contractible in $X$. By continuity of $u$ and compactness of $K$, there is $r_0>0$ such that $u(N_{r_0}(K))\subset N_{\varepsilon}(T)$. 

Assume that for any $0<r<r_0$ there exists a curve $\gamma_{r}$ in $N_{r}(K)$ which is not contractible in $M$. Any continuous map 1-homotopic to $\varphi$ is homotopic to $\varphi$ (see \cite[Lemma~6.2]{SW20}), as every surface not of sphere-type has trivial second homotopy group. This implies that the curve $u\circ\gamma_{r}\subset N_{\varepsilon}(T)$ is not contractible in $X$, a contradiction. Thus, we can choose $0<r<r_0$ such that any curve in $N_r(K)$ is contractible in $M$. We claim:

\bl
\label{lem:everyCurveContractible}
Let $K\subset M$ be a compact and connected set and let $r>0$ be such that every curve in $U:=N_r(K)$ is contractible in $M$. Then, the set $K$ is contained in the closure of a Jordan domain in $M$.
\el

With this lemma at hand we can easily finish the proof of the theorem. Indeed, from Lemma~\ref{lem:everyCurveContractible} we can deduce the existence of a Jordan domain $\Omega\subset M$ such that $K\subset\overline{\Omega}$. As $K$ is not cell-like, there exists a connected component of $\overline{\Omega}\setminus K$ that does not intersect $\partial\Omega$. In particular, there exists a connected component $U\subset\overline{\Omega}$ of $M\setminus u^{-1}(T)$ not intersecting $\partial M$. By arguing as in the proof of \cite[Lemma~4.2]{MW21}, we may assume that $\Omega$ is bounded by a biLipschitz curve and thus is a Lipschitz domain. Moreover, $T$ is an absolute Lipschitz retract as $T$ is a single point or biLipschitz homeomorphic to a closed interval. In particular, there exists a Lipschitz retraction $P\colon X\to T$. Define $w:=P\circ u$. By arguing as in the proof of \cite[Theorem~1.2]{LW20-param}, there exists a map $u_1\in N^{1,2}(\Omega,X)$ having the same trace as $u|_{\Omega}$ and agreeing with $u$ on $\Omega\setminus U$ and with $w$ on $U$. After applying the general gluing theorem for Sobolev maps, see \cite[Theorem~12.1.3]{KS93}, we obtain a map $u_2\in\Lambda(M,\partial X,X)$ agreeing with $u$ on $M\setminus U$ and with $w$ on $U$. As both $U$ and $T$ are contractible in their respective spaces, we obtain that $u_2$ is homotopic to $\varphi$. Thus, $u_2\in \Lambda(M,\varphi,X)$.

The set $T$ is a single point or biLipschitz homeomorphic to a closed interval implying that $\apmd w_z$ is degenerate for almost every $z$. Hence, the inscribed Riemannian area of $w|_U$ is zero. Since $u$ minimizes the inscribed Riemannian area among all maps in $\Lambda(M,\varphi,X)$, it follows that the inscribed Riemannian area of $u|_U$ is zero. By \cite[Proposition~1.1]{FW20}, the Reshetnyak energy of $u|_U$ is zero as well and it follows that $u|_U$ is constant. Therefore $u(U)$ is contained in $T$, a contradiction. Thus every connected component of $u^{-1}(T)$ is cell-like and $u$ satisfies (ii). This finishes the proof of Theorem~\ref{thm:cont-and-inf-iso-implies-monotone}.\\

We are left to prove the lemma above. For a topological space $Y$, a subset $A\subset Y$ and two distinct points $x,y\in Y\setminus A$, we say that $A$ \emph{separates} $x$ from $y$ if every connected subset $B\subset Y$ containing $x$ and $y$ intersects $A$. Note that Lemma~\ref{lem:everyCurveContractible} holds for all surfaces that are not of disc- or sphere-type. Hence, we provide a proof also for surfaces with empty boundary.

\begin{proof}[Proof of Lemma \ref{lem:everyCurveContractible}] We first assume that $M$ is a closed surface. Equip $M$ with a Riemannian metric $g$ and let $\phi\colon\widehat{M}\to M$ be the universal cover of $M$. Fix $x_0\in U$. The fiber of $U$ under $\phi$ is given by
$$\phi^{-1}(U)=\bigcup_{\alpha\in\pi_1(M,x_0)}U_{\alpha},$$
where we define 
$$U_{\alpha}:=\{[\alpha+\beta]:\beta\colon[0,1]\to U \text{ continuous, } \beta(0)=x_0\}.$$
As every curve in $U$ is contractible we have $U_{\alpha}\cap U_{\alpha'}=\emptyset$ for $\alpha\neq\alpha'$ and $ \phi|_{U_{\alpha}}\colon U_{\alpha}\to U$ is a homeomorphism for every $\alpha\in\pi_1(M,x_0)$. Fix $\alpha\in\pi_1(M,x_0)$. The set $U_{\alpha}$ is bounded. Indeed, if $U_{\alpha}$ is not bounded, there exists a curve $\tau$ in $U$ winding infinitely many times around a handle of $M$. Since every point in the image of $\tau$ is contained in a ball of radius $r$ in $U=N_r(K)$, we can find a non-contractible curve contained in $U$, a contradiction. Hence, there exists a set $Z\subset\widehat{M}$ homeomorphic to $\overline{D}$ with $U_{\alpha}\subset Z$. Observe that the set $K_{\alpha}:={ \phi|_{U_{\alpha}}}^{-1}(K)$ is homeomorphic to $K$. 
The coarea inequality for Lipschitz maps (see e.g. \cite[Theorem~2.10.25]{Fed69}) implies that the set 
$$S_t:=\{y\in Z:\dist(y,K_{\alpha})=t\}$$
has finite $\mathcal{H}^1$-measure for almost every $t>0$. Let $t\in(0,r)$ be such that $\mathcal{H}^1(S_t)<\infty$. The set $S_t$ separates any point $p\in K_{\alpha}$ from any $q\in\partial Z$ and, by \cite[Corollary~7.6]{LW16-intrinsic}, contains a Jordan curve $\gamma_t$ still separating $p$ from $q$.
By the Jordan curve theorem, $\gamma_t$ bounds a Jordan domain $\Omega_t$ in $Z$ with $K_{\alpha}\subset\Omega_t$. 

Consider the set $\gamma:= \phi(\gamma_t)\subset U$, which is homeomorphic to $\gamma_t$ and hence a Jordan curve. By assumption, every curve in $U$ is contractible in $M$ and therefore $\gamma$ bounds a Jordan domain $\Omega\subset M$. As $\overline{\Omega}$ is contractible, the map $ \phi|_{\overline{\Omega}_{\alpha}}\colon\overline{\Omega}_{\alpha}\to\overline{\Omega}$ is a homeomorphism, where $\overline{\Omega}_{\alpha}$ is defined analogously to $U_{\alpha}$. Hence, $\overline{\Omega}_{\alpha}$ is homeomorphic to $\overline{D}$ and shares the same boundary as $\overline{\Omega_t}$. Assume that $K$ is not contained in $\overline{\Omega}$. Since $ \phi|_{\overline{\Omega}_{\alpha}}$ is a homeomorphism and $\partial\overline{\Omega}_{\alpha}=\gamma_t$ does not intersect $K_{\alpha}$, it holds that $K_{\alpha}\subset \widehat{M}\setminus\Omega_{\alpha}$. This implies $\overline{\Omega}_{\alpha}\neq\overline{\Omega_t}$ and therefore $\overline{\Omega}_{\alpha}\cup\overline{\Omega_t}$ is a sphere, a contradiction.

If $M$ has non-empty boundary we consider the Schottky double $(M^*,g^*)$ of $(M,g)$, obtained by gluing two copies of $M$ along their boundaries and by doubling the metric $g$; compare to the proof of \cite[Lemma 2.4]{FW21}. By construction, $(M^*, g^*)$ is a smooth closed surface containing an isometric copy of $M$, denoted again by $M$. We use the same strategy as above to obtain the existence of a Jordan domain $\Omega$ in $M^*$ containing $K\subset M\subset M^*$. A connected component of the intersection of $\Omega$ with $M$ is again a Jordan domain whose closure contains $K$.
\end{proof}

\section{Applications}\label{sec:applications}
In this short section we briefly describe how a quasiconformal almost parametrization upgrades to a quasisymmetric map under the assumptions of Ahlfors $2$-regularity and linear local connectedness and to a geometrically quasiconformal map after assuming reciprocality. In particular, we show that Theorem~\ref{thm:main} implies generalizations of the uniformization theorems of Bonk and Kleiner as well as Rajala.

A homeomorphism $f\colon X\to Y$ between metric spaces is \emph{quasisymmetric} if there exists a homeomorphism $\eta:[0,\infty)\to[0,\infty)$ such that 
$$d_Y(f(x),f(y))\leq \eta(t)\cdot d_Y(f(x),f(z))$$
for all points $x,y,z\in X$ with $d_X(x,y)\leq t\cdot d_X(x,z)$. Moreover, a metric space $X$ is said to be \emph{Ahlfors $2$-regular} if there exists $K>0$ such that for all $x\in X$ and $0<r<\diam X$, we have
$$K^{-1}\cdot r^2\leq \hm^2(B(x,r))\leq K\cdot r^2.$$
We say that $X$ is \emph{linearly locally connected (LLC)} if there exists a constant $\lambda\geq1$ such that for all $x\in X$ and $r>0$, every pair of distinct points in $B(x,r)$ can be connected by a continuum in $B(x,\lambda r)$ and every pair of distinct points in $X\setminus B(x,r)$ can be connected by a continuum in $X\setminus B(x,r/\lambda)$.  

Note that every compact Ahlfors $2$-regular metric space is in particular of finite Hausdorff $2$-measure. Denote by $\Lambda(M,X)$ the family of Newton-Sobolev maps $u\in N^{1,2}(M,X)$ such that $u$ is a uniform limit of homeomorphisms from $M$ to $X$. Theorem~\ref{thm:main} shows that $\Lambda(M,X)$ is not empty for $M$ having non-empty boundary and $X$ being geodesic, Ahlfors $2$-regular and homeomorphic to $M$. By arguing exactly as in \cite[Section~5]{FM21}, there exists a canonical quasisymmetric homeomorphism from $M$ to $X$, if $X$ is furthermore LLC. After applying gluing techniques as in \cite{FM21}, we receive that this also holds for closed surfaces. Hence, we obtain the following theorem which recovers \cite[Theorem~1.1]{FM21} and generalizes Bonk-Kleiner's theorem \cite[Theorem~1.1]{BK02}.

\bt 
Let $X$ be a geodesic metric space which is Ahlfors $2$-regular, linearly locally connected and homeomorphic to a smooth surface $M$. 
Then, there exist a map $u\in \Lambda(M,X)$ and a Riemannian metric $g$ on $M$ such that
$$E_+^2(u,g)=\inf\{E_+^2(v,h):v\in \Lambda(M,X),\, h\text{ a smooth Riemannian metric on M}\}.$$
Any such $u$ is a quasisymmetric homeomorphism from $M$ to $X$ and the pair $(u,g)$ is uniquely determined up to a conformal diffeomorphism $(M, g)\to (M,h)$.
\et 

The assumption of $X$ being geodesic is not needed if $X$ is closed, since every closed, LLC and Ahlfors $2$-regular metric surface is geodesic up to a biLipschitz change of metric (see \cite[Theorem~B.6]{Sem96} and \cite[Lemma~2.5]{BK02}).

We now turn to the generalization of Rajala's uniformization theorem. Consider a metric space $X$ homeomorphic to a smooth surface $M$ and of finite Hausdorff 2-measure. A homeomorphism $u\colon M\to X$ is \emph{geometrically quasiconformal} if there exists $K\geq 1$ such that 
\begin{align}\label{eq:geom-qc-def-intro}
  K^{-1}\cdot \MOD(\Gamma)\leq \MOD(u\circ\Gamma)\leq K\cdot \MOD(\Gamma)
\end{align}
for every family $\Gamma$ of curves in $M$ with respect to a Riemannian metric $g$ on $M$. We call $X$ a \emph{quasiconformal surface} if every point of $X$ is contained in a geometrically quasiconformal image of $\overline{D}$. By Rajala's uniformization theorem \cite{Raj14}, this is equivalent to being \emph{locally recirocal}, i.e. every point of $X$ is contained in a reciprocal neighbourhood $U$ that is homeomorphic to $\overline{D}$. We call $U$ \emph{reciprocal} if the following two conditions hold. For every $x\in U$ and $R>0$ with $U\setminus B(x,R)\not=\emptyset$ we have 
\begin{align}\label{eq:points-zero-modulus}
 \lim_{r\to 0} \MOD(B(x,r), X\setminus B(x,R); \bar{B}(x,R)) = 0,
\end{align}
where $\MOD(E,F;G)$ denotes the modulus of the family of curves joining $E$ and $F$ in $G$ for some subsets $E,F,G\subset U$. Moreover, there exists $\kappa>0$ such that every closed topological square $Q\subset U$ with boundary edges $\zeta_1,\zeta_2,\zeta_3,\zeta_4$ in cyclic order satisfies
\begin{align}\label{eq:top-squares-moduli-opposite-edges}
 \MOD(\zeta_1, \zeta_3; Q)\cdot\MOD(\zeta_2,\zeta_4; Q)\leq \kappa.
\end{align}
Notice that Rajala \cite{Raj14} originally assumed an additional lower bound on the product in \eqref{eq:top-squares-moduli-opposite-edges}. It has been shown that this lower bound is always satisfied, see \cite{RR19} and \cite{EBPC21}. 

If $X$ is locally geodesic and has non-empty, rectifiable boundary, then by Theorem~\ref{thm:main}, there exists a continuous, monotone surjection $u\colon M\to X$ satisfying \eqref{ineq:modulus} with $K=\frac{4}{\pi}$. By compactness, there exists a constant $C\geq1$ and finitely many $C$-biLipschitz maps $\psi_i\colon\overline{D}\to M$ such that the sets $U_i:=\Image(u\circ\psi_i)$ are homeomorphic to $\overline{D}$ and cover $X$. If each $U_i$ satisfies \eqref{eq:points-zero-modulus}, then the maps $u\circ\psi_i\colon\overline{D}\to U_i$ are homeomorphisms (see \cite[Proposition~3.1]{MW21}) and, by \cite[Proposition~3.3]{MW21}, upgrade to geometrically quasiconformal maps if each $U_i$ satisfies \eqref{eq:top-squares-moduli-opposite-edges}. Since $\psi_i$ are biLipschitz and $U_i$ cover $X$, the map $u$ itself is geometrically quasiconformal. We have thus established the following version of \cite[Theorem~1.2]{Iko19} for locally geodesic surfaces with non-empty boundary.

\bt
	Every locally geodesic quasiconformal surface with $k\geq 1$ rectifiable boundary components is quasiconformally equivalent to a Riemannian surface wit $k$ boundary components.
\et
{\setstretch{1}
\def\cprime{$'$} \def\cprime{$'$} \def\cprime{$'$} \def\cprime{$'$}

}

\end{document}